\documentstyle[12pt]{article}
\newtheorem{vtheorem}{Theorem}
\newtheorem{theorem}{Theorem}[section]
\newtheorem{corollary}[theorem]{Corollary}
\newtheorem{example}[theorem]{Example}
\newtheorem{definition}[theorem]{Definition}
\newtheorem{proposition}[theorem]{Proposition}

\newtheorem{CD}[vtheorem]{Corollary/Definition}

\begin{document}

\title{On conditional expectations of finite index}
\author{Michael Frank and Eberhard Kirchberg}
\maketitle

\begin{abstract}
\noindent
For a conditional expectation $E$ on a (unital) C*-algebra $A$ there
exists a real number $K \geq 1$ such that the mapping $K \cdot E - {\rm id}_A$
is positive if and only if there exists a real number $L \geq 1$ such that
the mapping $L \cdot E - {\rm id}_A$ is completely positive,
among other equivalent conditions. The estimate
$(\min \, K) \leq (\min \, L) \leq (\min \, K) [\min \, K]$ is valid, where
$[\cdot]$ denotes the integer part of a real number. As a consequence the
notion of a ''conditional expectation of finite index'' is identified with
that class of conditional expectations, which extends and completes results of
M.~Pimsner, S.~Popa [27,28], M.~Baillet, Y.~Denizeau and J.-F.~Havet [6] and
Y.~Watatani [35] and others. Situations for which the index value and the
Jones' tower exist are described in the general setting. In particular, the
Jones' tower always exists in the W*-case and for ${\rm Ind}(E) \in E(A)$ in
the C*-case. Furthermore, normal conditional expectations of finite index
commute with the (abstract) projections of W*-algebras to their finite,
infinite, discrete and continuous type I, type ${\rm II}_1$, type
${\rm II}_\infty$ and type III parts, i.e.~they respect and preserve these
W*-decompositions in full.
\end{abstract}

\medskip \noindent
{\it  Primary Classification:} 46L99    \newline
{\it Secondary Classification:} 46L10, 46L85, 46H25, 47C15.

\newpage

\noindent
Treating conditional expectations on C*-algebras we follow the definition
of M.~Takesaki's and {\c{S}}.~Str{\u{a}}til{\u{a}}'s monographs
\cite{Take,Strat}: Conditional expectations $E$ on a C*-algebra $A$ are
projections of norm one of $A$ into a C*-subalgebra $B \subseteq A$ leaving
$B$ invariant. Immediate consequences are: $E$ is a $B$-bimodule map on $A$,
and the extension of faithful conditional expectations $E$ to the bidual
Banach space and W*-algebra $A^{**}$ yield normal faithful conditional
expectations onto the canonically normally embedded W*-subalgebra $B^{**}$
preserving the common identity. Consequently, for faithful conditional
expectations on C*-algebras $A$ the C*-algebra can always supposed to be
unital, and $E(1_A) = 1_B = 1_A$, cf.~\cite{Take} and
\cite[Lemma 4.1.4]{Bing-Ren} for details.

\smallskip \noindent
Studying the literature about conditional expectations of finite
index on C*-algebras we obtain two mainstreams of investigations:
One direction is concerned with normal conditional expectations of finite
index on W*-factors and several phenomena which occur in more special
settings of this case.
The other direction treats the algebraically characterizable case (see
Y.~Watatani \cite{Wata}, \cite{Khosh}) where the original C*-algebra $A$ is
a finitely generated C*-module over the (unital) image C*-algebra $B$.

\smallskip \noindent
A first definition for conditional expectations to be of finite index
was given by M.~Pimsner and S.~Popa \cite{PP,Po} in the context of W*-algebras
$M$ generalizing results of H.~Kosaki \cite{Kosaki} and V.~F.~R.~Jones
\cite{Jones}:  There has to exist a constant $K \geq 1$
such that $(K \cdot E - {\rm id}_M)$ is a positive mapping on $M$. However,
attempts to describe the more general situation of conditional expectations
on C*-algebras with arbitrary centers to be ''of finite index'' in some
sense(s) get into difficulties, as the gap of knowledge separating
Proposition 3.3 and Theorem 3.5 of the paper \cite{BDH} on W*-algebras 
by M.~Baillet, Y.~Denizeau and J.-F.~Havet shows. For non-trivial centers of
W*-algebras $M$ they obtained that even in the case of normal faithful
conditional expectations $E$ on $M$ the index value can be calculated only in
situations when there exists a number $L \geq 1$ such that the mapping
$(L\cdot E - {\rm id}_M)$ is completely positive, which seems to be more
since the difference $(\min \, L) = (\min \, K)^2$ at least appears,
cf.~Example \ref{ex1}. S.~Popa showed that for normal conditional expectations
$E: M \to N \subseteq M$ the conditions of \cite[Prop.~3.3, Th.~3.5]{BDH}
are equivalent, \cite[Th.~1.1.6, Rem.~1.1.7]{Po}. Y.~Watatani's attempt to
overcome this difficulty in the algebraic way considering the C*-algebra $A$
as a finitely generated projective module over the (unital) image C*-algebra
$B$ of the conditional expectation $E$ turned out to be unsuitable to give a
general solution of the key problems.

\smallskip \noindent
Let us give an example which is characteristic for situations for which
the algebraic approach does not work but which are, nevertheless, well-behaved
in some sense: Let $A$ be the C*-algebra ${\rm C}([-1,1])$ all continuous
functions on the interval $[-1,1]$. Consider the normal conditional
expectation $E$ on $A$ defined by
  \[
  E(f)(x) = \frac{f(x)+f(-x)}{2} \,\, {\rm for} \,\, x \in [-1,1].
  \]
The image C*-algebra $B$ can be identified with the set $\{ f \in A : 
f(x)=f(-x) \: {\rm for} \: x \in [-1,1] \}$ inheriting the C*-structure from
$A$. The Hilbert $B$-module $\{ A, E(\langle .,. \rangle_A) \}$ is not
self-dual since the bounded $B$-linear mapping 
\[
E(f_0^*f)(x) =  \frac{f_0^*(x)f(x)+f_0^*(-x)f(-x)}{2} \: {\rm for} 
\: x \in [-1,1] \, , \, f \in A \quad ,
\]
\[
\quad f_0(x) = \left\{ \begin{array}{cc} 1 & x \in (0,1] \\
                                          0 & x=0 \\
                                         -1 & x \in [-1,0)
                 \end{array} \right.
\]
maps $A$ into $B$, however $f_0$ does not belong to $A$. Consequently, the
Hilbert $B$-module $\{ A, E(\langle .,. \rangle_A) \}$ cannot be finitely
generated and projective, and $E$ is not of finite index in the sense of 
Y.~Watatani's algebraic definition. Nevertheless, $E$ has a very good 
property: The mapping $(2 \cdot E - {\rm id}_M)$ is obviously completely
positive. Beside this a careful analysis yields the $B$-reflexivity of
the Hilbert $B$-module $\{ A, E(\langle .,. \rangle_A) \}$.

\smallskip \noindent
Investigating such kind of conditional expectations in some more detail we
obtain that the existence of a number $K \geq 1$ such that the mapping
$(K \cdot E - {\rm id}_A)$ is positive should be sufficient to imply very
good properties of $E$, see \cite[Prop.~3.3]{BDH}, \cite{AS},
\cite[Th.~3.4, 3.5]{Jol} and \cite[Th.~1.1.6, Rem.~1.1.7]{Po}. We use the
denotations
\[
K(E) = \inf \{ K : (K \cdot E - {\rm id}_A) \; {\rm is} \; {\rm positive}
\; {\rm on} \; A \} \, ,
\]
\[
L(E) = \inf \{ L : (L \cdot E - {\rm id}_A) \; {\rm is} \; \; {\rm completely}
\; {\rm positive} \; {\rm on} \; A \} \, .
\]
If there does not exist any finite number $K$ or $L$ with the properties
striven for, then $K(E)=\infty$ or $L(E)=\infty$.
To close the gap of knowledge separating Proposition 3.3 and Theorem 3.5 of
\cite{BDH} for the general C*-case and to obtain the right general definition
for conditional expectations on C*-algebras to be of finite index we show the
following fact generalizing the partial results by P.~Jolissaint
\cite[Th.~3.4, 3.5]{Jol}, by E.~Andruchow and D.~Stojanoff
\cite[Prop.~2.1, Cor.~2.4]{AS2} and by S.~Popa
\cite[Th.~1.1.6, Rem.~1.1.7]{Po}:

\begin{vtheorem} \label{theorem}
  Let $A$, $B$ be C*-algebras, where $B$ is a C*-subalgebra of $A$.
  Let
  \linebreak[4]
  $E: A \rightarrow B \subseteq A$ be a conditional expectation with fixed
  point set $B$.
  Then the following three conditions on $E$ are equivalent:
    \newcounter{cou1}
    \begin{list}{(\roman{cou1})}{\usecounter{cou1}}
      \item $E$ is faithful and the (right) pre-Hilbert $B$-module
            $\{ A, E(\langle .,. \rangle_A) \}$ is complete with respect to
            the norm $\|E(\langle .,. \rangle_A)\|^{1/2}_B$, where
            $\langle a,b \rangle_A=a^*b$ for every $a,b \in A$.
      \item There exists a number $K \geq 1$ such that the mapping
            $(K \cdot E - {\rm id}_A)$ is positive.
      \item There exists a number $L \geq 1$ such that the mapping
            $(L \cdot E - {\rm id}_A)$ is completely positive.
    \end{list}
  If the second condition is valid for some $K$, then every number $L \geq K
  \cdot [K]$ realizes the third condition, where $[K]$ denotes the integer
  part of $K$. In some situations the equality $L(E)=K(E)^2 \in {\bf N}$ is
  valid.
\end{vtheorem}

\begin{CD} {\rm (cf.~\cite{Kosaki,PP,Po}) \newline
  Let $A$ and $B$ be C*-algebras, where $B$ is a C*-sub\-al\-ge\-bra of
  $A$. Let $E: A \rightarrow B \subseteq A$ be a conditional expectation with
  fixed point set $B$. If there exists a number $K \geq 1$ such that the
  mapping $(K \cdot E - {\rm id}_A)$ is positive, then $E$ {\it is of finite
  index}. The index value can be calculated in the center of either the
  bidual W*-algebra $A^{**}$ or its discrete part.   }
\end{CD}

\noindent
{\it Plan of the proof:} \newline
We start with the simple observation that the conditional expectation $E$ on
the C*-algebra $A$ can be continued to a normal conditional expectation
$E^{**}$ on the bidual W*-algebra $A^{**}$ of $A$ preserving property (ii)
above with the same minimal number $K(E^{**})=K(E) \geq 1$. Normality follows
from a general extension property of conditional expectations
\cite[Lemma 4.1.4]{Bing-Ren}, whereas the additional property can be derived
from monotonicity or from the up-down-up theorem of G.~K.~Pedersen
\cite[II.2.4.]{Pe}, for example.

\noindent
The next step is explained separately in Proposition \ref{prop2} below:
The projection of W*-algebras to their discrete (i.e.~atomic type I) part
commutes with normal faithful conditional expectations $E$ of finite index
if condition (ii) holds for $E$, i.e.~the discrete part can be
considered separately without loss of generality. This leads to a detailed
study of discrete W*-algebras and conditional expectations on them.
The theorem will be proven for this class, see Propositions \ref{prop1},
\ref{prop2}.

\noindent
To return to the general C*-case we make use of a theorem of Ch.~A.~Akemann
stating that the $*$-homomorphism of a C*-algebra $A$ into the discrete 
part of its bidual W*-algebra $A^{**}$ which arises as the composition of
the canonical embedding of $A$ into $A^{**}$ followed by the projection to
the discrete part of $A^{**}$ is an injective $*$-homomorphism,
\cite[p.~278]{Ake1} and \cite[p.~I]{Ake2}.
A conditional expectation $E$ on $A$ with property (ii) for some real number
$K \geq 1$ has an extension $E^{**}$ to $A^{**}$ and in particular, to the
discrete part of $A^{**}$ where $A$ is faithfully $*$-represented. $E^{**}$
restricted to this $*$-representation of $A$ recovers $E$ and $K$, and
condition (ii) holds for the restriction of $E^{**}$ to the discrete part of
$A^{**}$ for the same number $K \geq 1$.
Since Theorem \ref{theorem} is fulfilled for $E^{**}$ it holds for its
restriction $E$ to the faithfully represented original C*-algebra $A$ with
the same minimal numbers $K(E),L(E)$, and condition (iii) and the estimates
turn out to be valid.

\smallskip \noindent
The equivalence of the conditions (i) and (ii) follows from
\cite[Prop.~3.3]{BDH} for the W*-case. 
To give an independent argument we show how to derive
condition (ii) from condition (i). If $A$ is complete with respect to the
norm $\|E(\langle .,. \rangle_A)\|^{1/2}$, then there exists a number
$K$ such that the inequality $K \|E(x^*x)\| \geq \|x^*x\|$ holds for every
$x \in A$ by the general theory of Banach spaces. Set $x=a(\varepsilon
+E(a^*a))^{-1/2}$ with $a \in A$ and observe that
\[
(\varepsilon+E(a^*a))^{-1/2} \cdot E(a^*a) \cdot
                                  (\varepsilon+E(a^*a))^{-1/2} \leq 1_A \, .
\]
This implies the inequality $K \cdot 1_A \geq (\varepsilon+E(a^*a))^{-1/2}
a^*a(\varepsilon+E(a^*a))^{-1/2}$,
and multiplying with $(\varepsilon+E(a^*a))^{1/2}$ from both sides
we obtain $K(\varepsilon+E(a^*a)) \geq a^*a$ for every $\varepsilon >0$,
every $a \in A$. This yields condition (ii). The converse is obvious by
spectral theory.
$\: \bullet$

\medskip \noindent
In the next four sections we explain some details and consequences of our
investigations. The first section is concerned with the property of normal
conditional expectations $E$ of finite index to commute with the abstract
projection of W*-algebras to their discrete part, as well as with the proof
of Theorem \ref{theorem} for the discrete W*-case. In section two we show that
those mappings $E$ commute with the abstract projections of W*-algebras to
their finite, infinite, continuous type I, type ${\rm II}_1$, type
${\rm II}_\infty$ and type III parts, too. This extends results of S.~Sakai
\cite{Sakai} and J.~Tomiyama \cite{Tomiy} for general normal conditional
expectations. Section three is devoted to the investigation of the
index value of general conditional expectations of finite index and of
Jones' tower constructions. There exists an index value for $E: A \rightarrow B
\subseteq A$ inside the center of $A$ iff the discrete part of the index value
of the extended normal conditional expectation $E^{**}$ inside the center of
the bidual W*-algebra $A^{**}$ of $A$ belongs to the canonical embedding of
$A$ into the discrete part of $A^{**}$. For normal conditional expectations
$E:M \rightarrow N \subseteq M$ of finite index on W*-algebras $M$ the Jones'
tower always exists. The Jones' tower exists in the general C*-case if
${\rm Ind}(E)$ is contained in the center of $B$. The last section collects
some interpretations of the results obtained in terms of non-commutative
topology and some dimension estimates in the case of finite centers.


\section{The discrete W*-case}

\bigskip \noindent
We want to consider {\it discrete} W*-algebras, i.e.~W*-algebras for
which the supremum of all minimal projections contained equals their
identity.
First of all, we have to recall the structure of normal conditional
expectations on type I factors over separable and finite dimensional
Hilbert spaces the image of which is a (type I) subfactor. As a partial
case we describe normal states. Subsequently we will make use of the inner
structure of these mappings.

\noindent
\begin{example}    \label{ex2} {\rm
   Let $E$ be a normal conditional expectation on the set $M= {\rm B}(l_2)$
   of all bounded linear operators on the separable Hilbert
   space $l_2$ or on $M = {\rm M}_n({\bf C})$, respectively.
   Suppose, its image is a W*-subfactor $N \subseteq M$. Then the Hilbert
   space $l_2$ (or ${\bf C}^n$) can be decomposed as the tensor
   product of two Hilbert spaces $H_o$ and $K_o$, $l_2 = H_o \otimes K_o$,
   (${\bf C}^n = H_o \otimes K_o$, resp.) such that the representation of $M$
   as the W*-tensor product ${\rm B}(H_o) \otimes {\rm B}(K_o)$ allows the
   description of $E$ by the formula
     \[
     E(T \otimes S) = T \otimes {\rm Tr}(C \cdot S) {\rm id}_{K_o}
     \]
   for a unique trace class operator $C \in {\rm B}(K_o)_h^+$ with ${\rm
   Tr}(C)=1$ and for arbitrary $T \in {\rm B}(H_o)$, $S \in {\rm B}
   (K_o)$ on elementary tensors, (cf.~\cite[Lemma 3.3.1, 3.3.2]{K.Frank},
   \cite[Prop.~2.4]{Tsui}).
   Note, that $E$ is faithful if and only if zero is not an
   eigenvalue of the operator $C$. To satisfy the condition (ii) of Theorem
   \ref{theorem} the inequality $\|T\| \leq K \|E(T)\|$ has to be valid
   for a positive number $K \in {\bf N}$ (to be fixed) and for every operator
   $T \in M$. Denote by $P_n \in {\rm B}(K_o)$ that projection
   which maps $K_o$ onto the eigenspace of the $n$-th eigenvalue $\lambda_n$
   of the trace class operator $C$. Then this inequality can be rewritten as
   \begin{equation}  \label{for555}
     1 = \| {\rm id}_{H_o} \otimes P_n \| \leq K \| E({\rm id}_{H_o} \otimes
     P_n) \| = \lambda_n K \, .
   \end{equation}
   If $K_o$ is infinite dimensional, then the eigenvalues of $C$ form a
   sequence converging to zero. Such an assumption contradicts to the
   finiteness and universality of the constant $K$ in the inequality
   (\ref{for555}). Hence, $E$ satisfies condition (ii) of Theorem
   \ref{theorem} if and only if the dimension of $K_o$ is finite.
   \newline
   The image $N=E(M)$ can be identified with the W*-subalgebra
   ${\rm B}(H_o) \otimes {\rm id}_{K_o} \simeq {\rm B}(H_o)$. Therefore,
   the set
   \[
   \{ {\rm id}_{H_o} \otimes P_{n,i} : i=1,...,{\rm dim}(P_n(K_o))\, , \,
   n \in {\bf N} \} \cup
   \{ {\rm minimal} \; {\rm partial} \; {\rm isometries} \; {\rm between}
   \; {\rm them} \}
   \]
   is an orthogonal basis of the Hilbert $N$-module $\{ M, E(\langle .,.
   \rangle_M) \}$, where the projections $\{ P_{n,i}: i \in {\bf N} \}$
   vary over a set of pairwise orthogonal minimal subprojections of the
   corresponding eigenspace projection $P_n$ of $C$.
   In case we assume Theorem \ref{theorem},(ii) to be valid the index of $E$
   exists in the sense of Y.~Watatani \cite{Wata} and equals
   \[
   {\rm Ind}(E) = \left( \sum_{n=1}^{{\rm dim}(K_o)}
                  (\lambda_n)^{-1} \right) \cdot {\rm id}_M \, .
   \]
   }
\end{example}

\noindent
Note, that this special structure of some conditional expectations of finite
index on $\sigma$-finite type I factors is not preserved for similar ones over
non-separable Hilbert spaces, as we shall obtain during our investigations.
The structure of $E$ described above implies that ${\rm Ind}(E) \in \{1\} \cup
[4,\infty)$, whereas in the non-separable case the discrete series of
possible index values between one and four arises additionally.

\smallskip \noindent
Now we are going to investigate normal conditional expectations on W*-algebras
satisfying condition (ii) of Theorem \ref{theorem} and their restrictions
to the discrete part. We start with a proposition which was obtained
by J.~Tomiyama \cite{Tomiy} investigating normal conditional expectations.
For completeness, we include a short proof of it which is different from
J.~Tomiyama's:

\begin{proposition} \label{prop1} {\rm (J.~Tomiyama)} \newline
  Let $M$ be a discrete W*-algebra and $E:M \rightarrow N \subseteq M$
  be a normal faithful conditional expectation with respect to $N$.
  Then $N$ is a discrete W*-algebra.
\end{proposition}

\noindent
{\it Proof:} 
Let $P_N \in N$ be the projection which maps $N$ onto its discrete 
part. The projection $P_N$ is known to be contained in the center of
$N$. Since $E(1_M-P_N)$ is positive and belongs to the center of $N$, too,
it has a carrier projection $r \in N$ which is contained in the center of $N$.
Note that $(1_M-P_N)$ and $N$ commute inside $M$. Consider the W*-subalgebra
$(1_M-P_N)M(1_M-P_N)$ of $M$. Let $r_1$ be the smallest central projection
of $N$ satisfying the equality $nr_1(1_M-P_N)=n(1_M-P_N)$ for every $n \in N$.
Then $(1_M-P_N)=(1_M-P_N)r_1$ and $r_1 \leq r$. Applying $E$ we obtain
$E(1_M-P_N) = E(1_M-P_N)r_1$ and hence, $r_1=r$. That means, the
mapping $nr \rightarrow nr(1_M-P_N)$ is a faithful normal $*$-homomorphism
from $Nr$ into $(1_M-P_N)M(1_M-P_N)$.

\noindent
Furthermore, $E$ maps the discrete W*-algebra $(1_M-P_N)M(1_M-P_N)$ faithfully
onto $Nr$, and $Nr$ does not possess any discrete part. Fix a normal state
$f$ on $Nr$ with support projection $p_f$. The centralizer $c(f)$ of the
state $f$ inside $p_f(Nr)p_f$ coincides with the fixed point algebra of the
associated KMS-group $\sigma_f$ on $p_f(Nr)p_f$, and there exists a normal
faithful conditional expectation $E': p_f(Nr)p_f \to c(f)$ such that
$f \equiv f \circ E'$ on $p_f(Nr)p_f$, cf.~\cite{Take1}. 
Composing the normal conditional expectations $E$, $E'$ and the
normal state $f$ we obtain a normal state $f \circ E' \circ E = f \circ E$
on the discrete W*-algebra $p_f(1_M-P_N)M(1_M-P_N)p_f$ which possesses a
non-discrete centralizer since $c(f)$ is non-discrete and contained in it.
This is a contradiction, since the intersection of the centralizer of the
normal state $f \circ E$ with $p_f(1_M-P_N)M(1_M-P_N)p_f$ has to be
isomorphic to a direct integral of orthogonal direct sums of finite matrix
algebras over the discrete center of $p_f(1_M-P_N)M(1_M-P_N)p_f$, which can
be supposed to be finite by a special choice of $f$. (This structure of the
centralizers can be derived from the structure of normal states
on discrete W*-factors, cf.~Example \ref{ex2} and \cite{Take}.)
Consequently, $r=0$, $1_M=P_N=1_N$ inside $M$, and $N$ has to be discrete.
$\: \bullet$

\smallskip \noindent
We want to remark that we have always to be careful of the structure of
the centers of $M$ and of $N$ and of their interrelation. For example,
set $M = \sum_{k \in {\bf Z}} M_{2,(k)}({\bf C})$ and consider the
normal faithful conditional expectation
\begin{eqnarray*}
M = \sum_{k \in {\bf Z}} M_{2,(k)}({\bf C}) & \rightarrow &
N = \sum_{k \in {\bf Z}} {\bf C}_{(k)} \\
\left\{ \left( \begin{array}{cc}
        a_{(k)} & b_{(k)} \\
        c_{(k)} & d_{(k)} \end{array} \right)_{(k)} : k \in {\bf Z} \right\} &
\rightarrow &
\left\{ \lambda_{(k)} = \frac{d_{(k)}+a_{(k+1)}}{2} : k \in {\bf Z} \right\}
\, ,
\end{eqnarray*}
where the normal embedding of $N$ into $M$ is defined by the formula
\[
 \left\{ \left( \begin{array}{cc}  \lambda_{(k-1)} & 0 \\
                            0 & \lambda_{(k)} \end{array} \right)_{(k)} \, :
                            \, k \in {\bf Z} \right\}  \, .
\]
The centers of $M$ and of $N$ are both isomorphic to $l_\infty({\bf Z})$,
however their intersection is the smallest possible one -- the complex
multiples
of the identity. By the way, the index of this mapping equals $4 \cdot 1_M$.
The most unpleasant circumstance is the non-existence of a common central
direct integral decomposition of $M$ and $N$ commuting with the conditional
expectation $E:M \rightarrow N$ described above causing difficulties in
proving. As the referee pointed out to us the following proposition was
independently observed by S.~Popa \cite[1.1.2]{Po}:

\noindent
\begin{proposition} \label{prop2}
  Let $M$ be a W*-algebra and $E:M \rightarrow N \subseteq M$ be a normal
  conditional expectation leaving $N$ invariant for which there exists a
  number $K \geq 1$ such that the mapping $(K \cdot E - {\rm id}_M)$ is a
  positive mapping. Then the essential part of the preimage of the discrete
  part of $N$ is contained in the discrete part of $M$, and the image of the
  discrete part of $M$ is exactly the discrete part of $N$. That is, the
  projection to the discrete part of W*-algebras commutes with normal
  conditional expectations $E$ on them possessing this additional property.
\end{proposition}

\noindent
{\it Proof:} Recall, that $E$ is faithful by the additional
condition. Let $p \in N$ be a minimal projection of $N$ and let $q \leq p$
be a projection of $M$. Since $E$ is faithful $E(q) \not= 0$ and the
inequality $0 < E(q) \leq E(p) = p$ holds. Since $p$ is minimal inside $N$
there exists a number $\mu \in (0,1]$ such that $E(q)= \mu p$. That is,
  \[
  E(\mu^{-1} \cdot q)=p \, ,
  \]
and because of the additional condition on $E$ we obtain
  \[
  K \cdot p \geq K \cdot E(\mu^{-1} q) \geq \mu^{-1} q > 0
  \]
and the estimate $\mu \geq K^{-1}$.

\noindent
Suppose, $p \in N$ can be decomposed into a sum of pairwise orthogonal
(arbitrary) projections $\{ q_\alpha : \alpha \in I\}$ inside $M$.
Obviously, $q_\alpha \leq p$ for every $\alpha \in I$, and
  \[
  p = E(p) = \sum_{\alpha \in I} E(q_\alpha)
  = \left( \sum_{\alpha \in I} \mu_\alpha \right) p
  \geq \left( \sum_{\alpha \in I} K^{-1}_{(\alpha)} \right) p \, .
  \]
Consequently, the sum has to be finite and the maximal number of non-trivial
summands is $[K]$, the integer part of $K$. We see, the projections
$\{ q_\alpha : \alpha \in I\} \in M$ possess only finitely many
subprojections and hence, belong to the discrete part of $M$. This
shows that the discrete part of $N=E(M)$ is completely contained in
the discrete part of $M$.

\noindent
By the way we obtain that every minimal projection $p \in N$ has to be
represented only in a small part of the central direct integral decomposition
of $M$ with finite-dimensional center. Extending $p \in N \subseteq M$ by
the partial isometries of $N \subseteq M$ every W*-factor block of $N$ turns
out to be represented in a part of the central direct integral decomposition
of $M$ with finite-dimensional center. Conversely, $E$ maps each W*-factor
block of $M$ into a part of the central direct integral decomposition of $N$
with finite-dimensional center, cf.~\cite[Cor.~3.18]{BDH}.

\smallskip \noindent
Let $P_N$ denote the projection which maps $N$ onto its discrete 
part. Subject to the first part of the present proof $P_N$ belongs to the
discrete part of $M$ since $P_N$ is the supremum of all minimal
projections of $N$. Applying Proposition \ref{prop1} to the discrete part
$P_M \cdot M$ of $M$ (where $P_M$ denotes the (central) projection of $M$
which carries its discrete part) we obtain the equality $P_M=P_N$ in the
center of $M$ by the faithfulness of $E$.
$\: \bullet$

\noindent
\begin{corollary}  \label{cor99}
  Let $M$ be a discrete W*-algebra and $E:M \rightarrow N \subseteq M$
  be a normal conditional expectation for which there exists a
  number $K \geq 1$ such that the mapping $(K \cdot E - {\rm id}_M)$ is
  positive.
  Then:
    \begin{list}{(\roman{cou1})}{\usecounter{cou1}}
      \item The center of $M$ is finite-dimensional if and only if the
            center of $N$ is finite-dimensional.
      \item The positive part of the preimage of a minimal projection $p \in
            N$ is contained in the finite dimensional W*-subalgebra $pMp$ of
            $M$ generated by the minimal projections of $M$ which are
            subprojections of $p$. The dimension of $pMp$ is at most $[K]^2$,
            where $[K]$ denotes the integer part of $K$.
      \item For a minimal projection $p \in N$ the minimal projections 
            generating $pMp$ are mapped by $E$ to the set
            $\{ \mu \cdot p : K^{-1} \leq \mu \leq 1 \}$.
            The image of $pMp$ is ${\bf C}p$, i.~e.~$E$ acts on $pMp$ as
            a normal state.
      \item If $q_1,q_2$ are two orthogonal minimal projections of $M$, then
            either $E(q_1Mq_1) \equiv$ $\equiv E(q_2Mq_2)$ or $E(q_1Mq_1)
            \cap E(q_2Mq_2) = \{ 0 \}$.
    \end{list}
\end{corollary}

\noindent
We restrict our attention to the situation when $N$ is a type I W*-factor.
Then the center of $M$ has to be finite dimensional. If we consider
$M$ as a self-dual (right) Hilbert $N$-module with the $N$-valued inner product
$E(\langle .,. \rangle_M)$, (where $\langle a,b \rangle_M= a^*b$ for $a,b \in
M$), we can try to count the index of $E$ applying \cite[Thm.~3.5]{BDH}. We
have to find a suitable Hilbert $N$-module basis $\{m_\alpha : \alpha \in I \}$
of $M$ and to count the sum $\sum_\alpha m_\alpha m_\alpha^*$. If this sum is
finite in $M$ in the sense of w*-convergence for some basis, then it is the
same for every other Hilbert $N$-module basis of $M$.

\noindent
We start with a decomposition of the identity $1_N=1_M$ into a w*-sum of
pairwise orthogonal minimal projections $\{ p_\nu \} \subset N$, and with a
subdecomposition of this sum into a w*-sum of pairwise orthogonal minimal
projections $\{ q_\alpha : \alpha \in I \} \subset M$. In our special setting
a suitable basis contains this maximal set of pairwise orthogonal minimal
projections $\{ q_\alpha \}$ of
$M$ weighted down by the inverse of the number $\mu_\alpha$ arising by the
equality $E(q_\alpha) = \mu_\alpha p_\nu$ for some minimal projection $p_\nu$
of $N$ with $\mu_\alpha \in [K^{-1},1]$, (compare with Example \ref{ex2}):
\[
\{ \sqrt{\mu_\alpha^{-1}} \cdot q_\alpha : \alpha \in I \} \subseteq
{\rm basis} \, .
\]
If $M$ is commutative, then the Hilbert $N$-module basis
of $M$ is complete, and
\[
{\rm Ind}(E) = \sum_{\alpha \in I} \sqrt{\mu_\alpha^{-1}} \cdot q_\alpha
                     \cdot (\sqrt{\mu_\alpha^{-1}} \cdot q_\alpha)^*
            \leq K \sum_{\alpha \in I} q_\alpha
            \leq K \cdot 1_M
\]
by \cite[Thm.~3.5]{BDH}. But, if $M$ is non-commutative, then we have to add
all the minimal partial isometries $\{ u_\beta : \beta \in J \}$ of $M$ each
connecting two minimal projections, but also weighted down by the inverse of
the number $\mu_\beta$ arising by the equality $E(u_\beta^*u_\beta)=
\mu_\beta p_\nu$ for a minimal projection $p_\nu$ of $N$, $\mu_\beta \in
[K^{-1},1]$, (compare with Example \ref{ex2} again):
\[
\{ \sqrt{\mu_\alpha^{-1}} \cdot q_\alpha : \alpha \in I \} \cup
\{ \sqrt{\mu_\beta^{-1}} \cdot u_\beta  : \beta \in J \}
\equiv {\rm basis} \, .
\]
The Hilbert $N$-module basis of $M$ is complete, but rather big. However, based
on special Hilbert W*-module isomorphisms (\cite{Fr2}) we can reduce the
generating set of partial isometries $\{ u_\beta : \beta \in J \} \subset M$.
For this aim we define equivalence classes of them by the rule: $\, u_\beta
\sim u_\gamma$ if and only if $q_\alpha u_\beta=q_\alpha u_\gamma \not= 0$ for
some minimal projection $q_\alpha \in M$ of our choice and $u_\beta = u_\gamma
v$ for some partial isometry $v \in N$ linking the projections $\{ p_\nu \}
\subset N$.
By application of Hilbert W*-module isomorphisms we obtain a sufficiently
large set of generators of $M$ as right Hilbert $N$-module if we select only
one representative of each equivalence class supplementary to our choice of
minimal projections. As the result the number of selected partial isometries
$\{ u_\beta : \beta \in J \} \subset M$ satisfying $q_\alpha u_\beta = u_\beta$
for a fixed minimal projection $q_\alpha \in M$ is limited by $([K]-1)$.

\noindent
Now we can easily estimate the index value in norm to check the
w*-convergence of the appropriate series inside $M$:
\begin{eqnarray} \label{schaetz}
{\rm Ind}(E) & = & \sum_{\alpha \in I} \sqrt{\mu_\alpha^{-1}} \cdot q_\alpha
                     \cdot (\sqrt{\mu_\alpha^{-1}} \cdot q_\alpha)^*  +
               \sum_{\beta \in J} \sqrt{\mu_\beta^{-1}} \cdot u_\beta
                     \cdot (\sqrt{\mu_\beta^{-1}} \cdot u_\beta)^*  \\
                     \nonumber
            & \leq & K \left( \sum_{\alpha \in I} q_\alpha
                        + \sum_{\beta \in J} u_\beta u_\beta^* \right)\\
                               \nonumber
            & \leq & (K \cdot [K]) \sum_{\alpha \in I} q_\alpha \\ \nonumber
            & \leq & (K \cdot [K]) \cdot 1_M   \, .   \\ \nonumber
\end{eqnarray}
Consequently, $\| {\rm Ind}(E) \|_M \leq K \cdot [K]$, and $E$ is of finite
index in the sense of M.~Baillet, Y.~Denizeau and J.-F.~Havet
\cite[Thm.~3.5]{BDH}. Furthermore, the mapping $((K \cdot [K]) \cdot E -
{\rm id}_M)$ is completely positive by the same theorem.  $\: \bullet$

\begin{proposition}  \label{prop3}
  If $M$ is supposed to be a discrete W*-algebra, then Theorem 1 is valid.
\end{proposition}

\noindent
For a proof we have only to realize that the equality
$L(E) = \| {\rm Ind}(E) \|$ is valid by \cite[Th.~3.5, (a)+(b)]{BDH}, and
the value $\| {\rm Ind}(E) \|$ was estimated by the number $(K(E) \cdot
[K(E)])$ from above by (\ref{schaetz}). (Of course, in special situations
this estimate may be far from being sharp.)

\begin{example}  \label{ex1} {\rm \cite[Exemples 3.7]{BDH}
\newline
  Let $N$ be a W*-algebra. Let $M$ be the C*-algebra
  of all 2$\times$2-matrices with entries from $N$. The embedding of $N$
  into $M$ can be described as that subset of $M$ consisting of the
  $N$-multiples of the identity matrix. Consider the conditional expectation
     \[
     E \left( \begin{array}{cc} a&b\\c&d \end{array} \right) =
     \frac{a+d}{2} \left( \begin{array}{cc} 1&0\\0&1 \end{array} \right)
     \, .
     \]
  Denote by $e_{jk}$ those elements of $M$ which possess only one non-zero
  element at that place where the j-th row and the k-th column intersect, the
  identity of $N$. Then the set
  \linebreak[4]
  $\{ \sqrt{2} \cdot e_{jk} : j,k=1,2 \}$ is a Hilbert $N$-module basis of $M$,
  and the index of $E$ has the value ${\rm Ind}(E)=4 \cdot {\rm id}_{M_2}$.
  \newline
  If $N$ is commutative, then the mapping $(K(i \circ E) -{\rm id}_M)$ is 
  positive for $K \geq 2$ already, whereas it is completely positive for
  $L \geq 4$ only. That is, the minimal constants of item (ii) and of item
  (iii) of Theorem \ref{theorem} may be different, in general, and in our
  special setting $4 = L(E) = K(E)^2 = 2^2 \in {\bf N}$.
  \newline
  At the contrary,
  if $N$ is a type ${\rm I}_\infty$, type ${\rm II}_1$ or separably
  representable infinite W*-factor, then $L(E) = K(E)$, see \cite{BDH,Wata}.
  That is, our estimate of the number $L(E)$ by the number $K(E)$ does not
  give a formula to calculate the value of $L(E)$ precisely, in general.

  \noindent
  Furthermore, a general estimate $L(E) \leq [K(E)]^2$ is not true: Modify
  the preceding example with $\lambda \in (0,1)$ setting
     \[
     E_\lambda \left( \begin{array}{cc} a&b\\c&d \end{array} \right) =
     ( \lambda a + (1-\lambda) d )
     \left( \begin{array}{cc} 1&0\\0&1 \end{array} \right)
     \, .
     \]
  Then $K(E_\lambda) = \max \{ \lambda^{-1}, (1-\lambda)^{-1} \}$ and
  ${\rm Ind}(E_\lambda) = L(E_\lambda) = \lambda^{-1} + (1-\lambda)^{-1}$.
  For every small $\varepsilon > 0$ and the choice $\lambda(\varepsilon) =
  1/(2 + \varepsilon)$ the assumption $L(E_\lambda) \leq
  [K(E_\lambda)]^2$ leads to the contradiction $\varepsilon^2 < 0$ since
  \[
  L(E_{\lambda(\varepsilon)}) = 4 + \frac{\varepsilon^2}{1+\varepsilon} \: ,
  \: K(E_{\lambda(\varepsilon)}) = 2 + \varepsilon \, .
  \]
  }
\end{example}

\noindent
\begin{corollary}
In contrast to the separable situation (see Example \ref{ex2})
there exist type I factors $M$ on non-separable Hilbert spaces and normal
conditional expectations $E$ of finite index on them with factor image
realizing the discrete index series between one and four.
\end{corollary}

\noindent
{\it Proof:} Start with the hyperfinite type ${\rm II}_1$ factor $A$ and the
appropriate (normal) conditional expectation $E$ of finite index on it
the image W*-algebra of which is a factor again, (see \cite{Jones,PP}).
Turn to the discrete part $M$ of the bidual W*-algebra $A^{**}$ of $A$.
The extended normal conditional expectation $E^{**}$ has the same invariant
$K(E^{**})=K(E)$, hence, also its restriction to the discrete part $M$ of
$A^{**}$ does. But $M$ is a type I factor, i.e.~$M$ has not any non-trivial
central projection since
$A$ has not any non-trivial two-sided norm-closed ideal. Indeed, every central
projection $p \in M$ can be decomposed into a orthogonal sum of minimal
projections $\{ p_\alpha \}$. For every projection $p_\alpha$ the orthogonal
complement $(1_M-p_\alpha)$ is the carrier projection of the w*-closure of
an appropriate maximal norm-closed left (or right) ideal $I_\alpha$ of $A$.
Consequently, $(1_M-p) \in {\rm Z}(M)$ should be the carrier projection of the
intersection $\cap_\alpha I_\alpha \equiv \cap_\alpha I_\alpha^* \subseteq A$
which should be a norm-closed two-sided ideal, and the claim is obvious.
$\: \bullet$

\medskip \noindent
As a corollary we obtain the structure of the relative commutant in the general
C*-setting:

\begin{corollary}  \label{cor89}
  Let $E:A \rightarrow B \subseteq A$ be a conditional expectation of finite
  index on a C*-algebra $A$. The relative commutant of $B$ inside $A$ is a
  subhomogeneous C*-algebra of finite type, i.e.~it is a C*-subalgebra
  of some matrix algebra ${\rm M}_n(C)$ with $n < \infty$ and with entries
  from a commutative W*-algebra $C$, for example $C={\rm Z}
  ((B^{**})_{discr.})$.
\end{corollary}

\noindent
{\it Proof:}
Consider the restriction of the (normal) extended conditional expectation
$E^{**}$ of $E$ to the discrete part of $A^{**}$. Since $E$ maps the relative
commutant $B' \cap A$ of $B$ with respect to $A$ to the center of $B$
the same is true for the restricted mapping $E^{**}$ and the appropriate
injectively $*$-represented C*-subalgebras of the discrete part of $A^{**}$.
Note, that $(B^{**})' \cap A^{**} \supseteq B' \cap A$. The center of the
discrete part of $B^{**}$ is discrete, and the preimage of every minimal
projection of it is a matrix algebra of dimension lower-equal $[K(E)]^2$ by
Corollary \ref{cor99},(i). Consequently, the injectively $*$-represented
C*-algebra $B' \cap A$ is contained in the C*-algebra of all
$n \times n$-matrices with $n = [K(E)]^2$ and with entries from the
center of the discrete part of $B^{**}$.
$\: \bullet$


\section{Decomposition preserving properties}

\medskip \noindent
By the work of S.~Sakai \cite{Sakai} and of J.~Tomiyama \cite{Tomiy} it has
been known that normal conditional expectations preserve
semi-finiteness, type I and discreteness of W*-algebras as properties of
their image W*-subalgebras. We want to show that normal conditional
expectations $E$ of finite index not only commute with the abstract projection
of W*-algebras to their discrete part but also with many canonical abstract
projections of W*-algebras to other parts of them. We obtain that such mappings
$E$ possess a decomposition into their restrictions to the appropriate type
components of the W*-algebras involved. The same observations were made by
S.~Popa \cite[1.1.2]{Po} independently, a fact which was brought to our
attention by the referee. Another result on projections which belong to
both ${\rm Z}(N)$ and ${\rm Z}(M)$ can be found at \cite[Th.~1]{FiIs}.

\noindent
We would like to point out that the word 'preimage' subsequently refers to
that part of the preimage of $N$ via $E$ which is spanned by the positive
part of it inside $M$. Of course, the kernel of $E$ is spread out over all
parts of the central decomposition of $M$ into W*-types, in general.

\noindent
\begin{proposition} \label{prop4}
  Let $M$ be a W*-algebra and $E:M \rightarrow N \subseteq M$ be a normal
  faithful conditional expectation leaving $N$ invariant. Then the image of
  the finite part of $M$ is exactly the finite part of $N$. 
  \newline
  If $E$ is additionally of finite index, then the preimage of the finite
  (resp., infinite) part of $N$ is contained in the finite (resp., infinite)
  part of $M$. That is, the projection to the finite (resp., infinite) part
  of W*-algebras commutes with normal conditional expectations $E$ of finite
  index. 
  Moreover, $E$ commutes with the projection to the type ${\rm II}_1$ part,
  to the continuous type ${\rm I}_{fin.}$ part and with the projection to the
  non-discrete infinite part of W*-algebras.
\end{proposition}

\noindent
{\it Proof:}
Suppose $E$ to be faithful. Denote a normal faithful center-valued
normalized trace on the finite part $M \cdot P_{fin}$ of $M$ by ${\rm tr}(.)$.
Suppose, $r \in N$ is the central carrier projection of $E(P_{fin})$.
By the same arguments as in the proof of Proposition \ref{prop1} the
C*-subalgebras $Nr$ and $N\cdot E(P_{fin})$ are identical W*-subalgebras.
The mapping $E \circ {\rm tr} \circ E$ is a normal faithful positive
normalized mapping on $Nr$ taking values in the center of $Nr$. Beside this,
$E \circ {\rm tr} \circ E$ is tracial on $Nr$. Commutative W*-algebras are
generally known to possess a decomposition into $\sigma$-finite W*-subalgebras
like
\[
\sum_{\alpha \in I} \oplus {\bf C}_{(\alpha)} \: \: \oplus  \: \:
         \sum_{\beta \in J} \oplus L^\infty ([0,1], \lambda)_{(\beta)}
\]
up to $*$-isomorphism by \cite[III.1.22]{Take} and
\cite[Prop.~1.14.10]{Bing-Ren}, where $\lambda$ denotes the Lebesgue measure
on $[0,1]$.
Considering the direct integral decomposition of $Nr$ over its center and
taking faithful normal states on the $\sigma$-finite components of
${\rm Z}(Nr)$ to be composed with the mapping $E \circ {\rm tr} \circ E$ we
obtain a separating set of tracial normal states on $Nr$. Consequently, $Nr$
is finite, cf.~\cite[Prop.~6.5.15]{Bing-Ren}.

\smallskip \noindent
The second statement is an easy consequence of the proof of
\cite[Th.~1.6,(1)]{Jol} if we consider arbitrary finite W*-algebras $N$
and apply the criterion on finiteness of W*-algebras of Sakai
(\cite[Th.~6.3.12]{Bing-Ren}) in this generalized setting.

\smallskip \noindent
The continuous type ${\rm I}_{fin.}$ part of $M$ (i.e.~the part with
continuous center and finite-dimensional fibres in the appropriate direct
integral decomposition over it) is mapped to the continuous type
${\rm I}_{fin.}$ part of $N$ since $E$ maps the type ${\rm I}$ part of $M$ to
the type ${\rm I}$ part of $N$ by the results of J.~Tomiyama \cite{Tomiy} and
since $E$ preserves the finite part and the discrete part by our previous
results. Conversely,
the projection of $N$ to its continuous type ${\rm I}_{fin.}$ part can be
decomposed into the sum of finitely many pairwise orthogonal abelian
projections $\{ p_n : n=1,...,m \} \subseteq N$. Since
\[
E(x)= \sum_{n_1,n_2 =1}^m p_{n_1}E(x)p_{n_2} =
\sum_{n_1,n_2 =1}^m E(p_{n_1}xp_{n_2})
\]
for all $E(x)$ of the continuous type ${\rm I}_{fin.}$ part of $N$, and
since the W*-subalgebras $\{ p_{n_1}Np_{n_2} : n_1,n_2=1,...,m \}$ are
commutative the preimage of the continuous type ${\rm I}_{fin.}$ part of $N$
has to be of finite subhomogeneous type inside $M$ by Corollary \ref{cor89}.
That is, it is contained in the continuous type ${\rm I}_{fin.}$ part of $M$
since $E$ preserves the discrete and the finite W*-parts.

\smallskip \noindent
Since the projections of W*-algebras to their discrete, their finite and their
infinite parts, respectively, commute with conditional expectations $E$
of finite index on them and because all they are realized by central
projections $P_{discr.}$, $P_{I_{cont.fin.}}$, $P_{fin.}$ and $P_{infin.}$ of
${\rm Z}(M) \cap {\rm Z}(N)$ the mapping $E$ commutes with the product
projections $P_{II_1} = P_{fin.} (1_M-P_{discr.}-P_{I_{cont.fin.}})$ and
$P_{infin.} (1_M-P_{discr.})$, too.
$\: \bullet$

\begin{proposition} \label{prop5}
  Let $M$ be a W*-algebra and $E:M \rightarrow N \subseteq M$ be a normal
  faithful conditional expectation leaving $N$ invariant. If $E$ is of finite
  index, then:

    \begin{list}{(\roman{cou1})}{\usecounter{cou1}}
      \item The image of the continuous type ${\rm I}_\infty$ part of $M$ is
            contained in the continuous type ${\rm I}_\infty$ part of $N$, and
            the preimage of the continuous type ${\rm I}_\infty$ part of $N$
            is contained in the continuous type ${\rm I}_\infty$ part of $M$.
      \item The image of the type ${\rm II}_\infty$ part of $M$ is contained
            in the type ${\rm II}_\infty$ part of $N$, and the preimage of
            the type ${\rm II}_\infty$ part of $N$ is contained in the
            type ${\rm II}_\infty$ part of $M$.
      \item The image of the type ${\rm III}$ part of $M$ is contained
            in the type ${\rm III}$ part of $N$, and the preimage of
            the type ${\rm III}$ part of $N$ is contained in the
            type ${\rm III}$ part of $M$.
    \end{list}
\end{proposition}

\noindent
{\it Proof:}  By \cite[Lemma V.2.29]{Take} the preimage of the type ${\rm III}$
part of $N$ is completely contained in the type ${\rm III}$ part of $M$
even if the faithful conditional expectation $E$ is not assumed to be of
finite index. This implies that the semi-finite part of $M$ is mapped by $E$
into the semi-finite part of $N$.

\smallskip \noindent
Let $E$ be of finite index. The type ${\rm II}_\infty$ part of $M$ is denoted
by $M \cdot P_{{\rm II}_\infty}$. There exists a faithful semi-finite normal
center-valued trace ${\rm tr}$ on $M \cdot P_{{\rm II}_\infty}$,
cf.~\cite[Prop.~6.5.8]{Bing-Ren}. Since $E$ commutes with the projection of
W*-algebras to their finite part and to their discrete part by Proposition
\ref{prop2} and \ref{prop4} the W*-subalgebra $E(M \cdot P_{{\rm II}_\infty})$
is contained in the direct sum of the continuous type ${\rm I}_\infty$ and of
the type ${\rm II}_\infty$ part of $N$. Consider the mappings
\[
\{ f \circ E \circ {\rm tr} \circ E \, : \: f \in {\rm Z}(N)_* \}
\]
on the W*-subalgebra $E(M \cdot P_{{\rm II}_\infty}) \subseteq N$. These
mappings form a faithful family of semi-finite normal traces on the positive
cone of $E(M \cdot P_{{\rm II}_\infty})$. By \cite[Prop.~6.5.7]{Bing-Ren}
the W*-algebra $E(M \cdot P_{{\rm II}_\infty})$ has to be semi-finite, and
hence, of type ${\rm II}_\infty$.

\smallskip \noindent
Conversely, let $N$ be of type ${\rm II}_\infty$ and $M$ arbitrary.
Let $p \in N$ be a finite projection. Then the preimage of the type
${\rm II}_1$ W*-subalgebra $pNp$ is exactly the W*-subalgebra $pMp$.
By Proposition \ref{prop4} $pMp$ has to be of type ${\rm II}_1$, too,
since $E$ is supposed to be of finite index.
By \cite[Th.~6.5.10]{Bing-Ren} there exists an increasing net $\{ p_\alpha :
\alpha \in I\}$ of finite projections in $N$ with strong limit $1_N$
inside $N$. Since $E$ is faithful and normal the least upper bound
of the net $\{ p_\alpha \} \subset N \subseteq M$ with respect to $M$ is
$1_M=1_N$, again.
Consequently, $M$ equals the w*-closure of the union of all type
${\rm II}_1$ W*-subalgebras $\{ p_\alpha M p_\alpha : \alpha \in I \}$,
and it has to be of type ${\rm II}_\infty$.

\smallskip \noindent
Now let $N$ be of continuous type ${\rm I}_\infty$ and $M$ arbitrary.
Let $p \in N$ be an abelian projection. The preimage of the commutative
W*-subalgebra $pNp$ of $N$ equals the W*-subalgebra $pMp \subseteq M$.
It has to be of subhomogeneous type by Corollary \ref{cor89},
and hence, it has to be contained in a matrix algebra of finite size with
entries from some commutative W*-algebra. That is, $pMp$ is of continuous
type ${\rm I}_{fin.}$ and possesses sufficiently many abelian subprojections
of $p$ with least upper bound $p$.
Since the set of all abelian projections $\{ p_\alpha : \alpha \in I\}$ of
$N$ has least upper bound $1_N$ with respect to $N$ the least upper bound of
them with respect to $M$ has to be $1_M=1_N$, too, since $E$ is normal and
faithful. Then $M$ is the w*-closure of the union of continuous type
${\rm I}_{fin.}$ W*-subalgebras $\{ p_\alpha M p_\alpha : \alpha \in I \}$.
The set of abelian projections of $M$ has least upper bound $1_M$ 
and hence, $M$ is of continuous type ${\rm I}_\infty$.

\smallskip \noindent
Summing up the considerations above $E$ maps the type ${\rm III}$ part of $M$
into the type ${\rm III}$ part of $N$, and the continuous type ${\rm I}_\infty$
part of $M$ into the continuous type ${\rm I}_\infty$ part of $N$.
$\: \bullet$

\medskip \noindent
Since ${\rm End}_N(M)$ and $N$ have a common center and since the W*-type
of $p \cdot {\rm End}_N(M)$ is the same as fore $pN$ for any suitable
$p \in {\rm Z}(N)$ we conclude that for finite W*-algebras $M$ with
finite-dimensional center and normal conditional expectations of finite
index $E: M \to N \subseteq M$ the W*-algebra $M$ is always a finitely
generated projective $N$-module in the standard decomposition of $M$ as
a Hilbert $N$-module via $E$, cf.~\cite[Th.~2.2,(2)]{Jol}.

\noindent
By the examples of P.~H.~Loi \cite{Loi} we conclude that conditional
expectations of finite index on type ${\rm III}_\lambda$ factors, 
$0 < \lambda < 1$, can map them to type ${\rm III}_{\lambda^m}$ subfactors
for any natural number $m \in {\bf N}$. Consequently, we can expect 
type preserving properties for conditional expectations of finite index
on type III W*-algebras only in the type ${\rm III}_0$ and ${\rm III}_1$
cases, if at all. This question has to remain unsolved for the time being.

\section{The general Jones' tower construction}

\noindent
The next step should be the construction of the Jones' tower in the general
C*-case and the estimate of a hopefully existing index value for conditional
expectations
\linebreak[4]
$E: A \rightarrow B \subseteq A$. The general construction of the Jones' tower
in the W*-case was shown by S.~Popa \cite[1.2.2]{Po}.
We are going to obtain a general index notion and a general Jones' tower
construction in some specific C*-cases. Unfortunately, the index value belongs
to the discrete part of $A^{**}$ only, in general, and not to the original
C*-algebra $A$ itself.

\smallskip \noindent
We consider the conditional expectation $E: A \rightarrow B \subseteq A$ of
finite index. The C*-algebra $A$ has the structure of a (right) Hilbert
$B$-module setting $\{A, E(\langle .,. \rangle_A) \}$, where
$\langle a,b \rangle_A =a^*b$ for $a,b \in A$ by Theorem \ref{theorem}.
The conditional expectation $E$ acts as a bounded (right-)B-linear mapping
on $A$. Furthermore, it can be identified with an elementary ''compact''
operator $e=\theta_{1_A,1_A}$ on the Hilbert $B$-module $A$ by the formula
$e(x)=\theta_{1_A,1_A}(x)=1_A \cdot E(1_A^*x) (=E(x))$. The projection $e$ is
the first projection to build the Jones' tower. It always exists since
the Hilbert $B$-module $B$ is a direct summand of the Hilbert $B$-module
$\{ A, E(\langle .,. \rangle_A ) \}$.

\noindent
The set of ''compact'' $B$-linear operators ${\rm K}_B(A)$ on the Hilbert
$B$-module $A$ is defined as the norm closure of the linear hull of the
elementary ''compact'' operators
\[
\{ \theta_{a_1,a_2} : \theta_{a_1,a_2}(x)=a_2 \cdot E(a_1^*x) \, ,
\: (x \in A) \} \, .
\]
There is also a (faithful) $*$-representation $\pi_E$ of $A$ in the set
of all bounded adjointable module operators ${\rm End}_B^*(A)$ on $A$ by
multiplication operators from the left
\[
\{ \pi_E(a) : \pi_E(a)(x)=ax \, , \: (x \in A) \} \, .
\]
Note, that not every bounded $B$-linear operator on the Hilbert $B$-module
$\{A, E(\langle .,. \rangle_A) \}$ should possess an adjoint bounded
$B$-linear operator on it, cf.~\cite[Th.~5.6,6.8]{Fr2} for criteria.
The unital C*-algebra ${\rm End}_B^*(A)$ of all bounded adjointable $B$-linear
operators on the Hilbert $B$-module $A$ is the multiplier C*-algebra of the
C*-algebra ${\rm K}_B(A)$ by \cite[Lemma 16]{Green} and \cite[Th.~1]{Kasp}.
Its center can be identified with the center of $B$ by the formula
\[
b \in {\rm Z}(B) \rightarrow \pi_E'(b) \in {\rm Z}({\rm End}_B^*(A)) \:\, ,
\: \pi_E'(b)(x)=xb \; , \, (x \in A)
\, .
\]
Obviously, $\pi_E(a_2) \theta_{1_A,1_A} \pi_E(a_1^*) = \theta_{a_1,a_2}$
for $a_1,a_2 \in A$, i.~e.~the emphasized operator $e$ gene\-rates all
elementary ''compact'' operators on $A$ by two-sided ideal operations with
respect to the $*$-represented C*-algebra $A$ inside ${\rm End}_B^*(A)$. The
linear hull of the elementary ''compact'' module
operators on $A$ coincide with the C*-algebra ${\rm K}_B(A)$, and the latter
coincides with its multiplier C*-algebra ${\rm End}_B^*(A)$) if and only if
$A$ is a projective finitely generated $B$-module, which recovers the algebraic
case described by Y.~Watatani, cf.~\cite[Prop.~1.1]{Fr2}, \cite{Wata,Fr}.

\smallskip \noindent
What about a conditional expectation $E_1: {\rm End}_B^*(A) \rightarrow
\pi_E(A)$ ? Unfortunately, we can only obtain a finite faithful operator-valued
weight
\[
F_1 : {\rm K}_B(A) \rightarrow \pi_E(A) \subseteq {\rm End}_B^*(A) \: \,  ,
\: F_1(\theta_{a_1,a_2}) = \pi_E(a_1a_2)  \, .
\]
In the algebraically characterizable case $F_1$ maps ${\rm id}_A \in
{\rm K}_B(A)$ to the index value ${\rm Ind}(E)$ of $E$ existing and
belonging to the center of $A$. The index value is greater-equal the
identity and hence, invertible. The sought for conditional expectation
$E_1$ arises as
\linebreak[4]
$E_1= \pi_E({\rm Ind}(E)^{-1}) \cdot F_1$ in that case.

\noindent
But, this easy construction does not work in the general case if the
identity operator on the Hilbert $B$-module $A$ is non-''compact''.
P.~Jolissaint gave a criterion in the W*-case showing that the right
construction can be obtained if and only if the number $L(E)$ is finite and,
hence, the index value exists for the normal conditional expectation $E$,
\cite[Prop.~1.5,(3)]{Jol}. By the results of S.~Popa \cite{Po} or by Theorem
\ref{theorem} this gives the general solution of the problem in the W*-case.
To overcome this difficulty in the general C*-case we make the following
definition for the index value of conditional expectations of finite index
and show its correctness afterwards:

\noindent
\begin{definition} {\rm
  Let $A$ be a C*-algebra and $E:A \rightarrow B \subseteq A$ be a 
  conditional expectation of finite index leaving $B$ invariant. The index
  value ${\rm Ind}(E)$ of $E$ is the projection of the index value
  ${\rm Ind}(E^{**})$ of the extended conditional expectation $E^{**}:
  A^{**} \rightarrow B^{**} \subseteq A^{**}$ to the discrete part of
  $A^{**}$.         }
\end{definition}

\begin{theorem}
  Let $A$ be a C*-algebra and $E:A \rightarrow B \subseteq A$ be a 
  conditional expectation of finite index leaving $B$ invariant.
  Then the index value ${\rm Ind}(E)$ is contained in the center of $A$,
  and there exists a conditional expectation $E_1: {\rm End}_B^*(A)
  \rightarrow \pi_E(A)$ mapping $e$ to $\pi_E({\rm Ind}(E)^{-1})$ if and only
  if ${\rm Ind}(E)$ belongs to the standardly embedded image of
  $A$ inside the discrete part of $A^{**}$.
  \newline
  For normal conditional expectations $E:M \rightarrow N \subseteq M$
  of finite index on W*-algebras $M$ the Jones' tower always exists.
  \newline
  The Jones' tower exists in the general C*-case if ${\rm Ind}(E)$ is
  contained in the center of $B$.
\end{theorem}

\noindent
{\it Proof:}
First, consider normal conditional expectations $E:M \rightarrow N \subseteq M$
of finite index on W*-algebras $M$. Then the normal conditional expectation
$E_1$ mapping the W*-algebra ${\rm End}_N(M)$ to $\pi_E(M)$ is faithful, and
there exists a number $K(E_1)=\| {\rm Ind}(E)^{-1} \|$ such that the mapping
$(K(E_1) \cdot E_1 - {\rm id}_{{\rm End}})$ is positive on ${\rm End}_N(M)$,
see \cite[Th.~3.5]{BDH}. But, by Theorem \ref{theorem} $E_1$ turns out to be
of finite index, and the Jones' tower can be built up repeating the basic
construction countably many times. This shows that the Jones' tower
construction always exists in the W*-case.

\smallskip \noindent
Now consider the general case. Suppose, ${\rm Ind}(E)$ belongs to $A$.
The finite faithful operator-valued weight
$F_1$ on the C*-algebra ${\rm K}_B(A)$ extends uniquely to a normal finite
faithful operator-valued weight $F_{1,{**}}$ on the C*-algebra
${\rm K}_{B^{**}}
(A^{**})$ by the way in which $E^{**}$ is derived from $E$. Obviously,
${\rm K}_B(A)$ can be considered as a C*-subalgebra of ${\rm K}_{B^{**}}
(A^{**})$, and $F_{1,{**}}$ restricted to ${\rm K}_B(A)$ recovers $F_1$.
If the projection of ${\rm Ind}(E^{**})$ to the discrete part of
$A^{**}$ is contained in the standard injective $*$-representation of $A$ in
the discrete part of $A^{**}$, then the restriction of $E_{1,{**}}$
(which exists as a normal conditional expectation from the W*-algebra
${\rm End}_{B^{**}}(A^{**})$ to $\pi_{E^{**}}(A^{**})$) to the C*-subalgebra
${\rm End}_B^*(A)$ of ${\rm End}_{B^{**}}(A^{**})$ gives a conditional
expectation $E_1$ from ${\rm End}_B^*(A)$ to $\pi_E(A)$ which equals $F_1$
multiplied by the inverse of the projection of ${\rm Ind}(E^{**})$ to the
discrete part of $A^{**}$, which belongs to ${\rm Z}(A)$ by assumption.

\noindent
Conversely, the extension $F_1^{**}$ of $F_1$ to the bidual linear space
and W*-algebra ${\rm K}_B(A)^{**}$ yields an image $F_1^{**}({\rm id})$ of the
identity inside the center of $\pi_E(A^{**})$. If there exists an extension
of $F_1$ to ${\rm End}_B^*(A)$ at all, then the projection of $F_1^{**}
({\rm id}) \in {\rm Z}(A^{**})$ to the discrete part of $A^{**}$ should equal
$F_1({\rm id}_A)$ because of the canonical embedding of multiplier C*-algebras
into bidual W*-algebras. Obviously, $F_1({\rm id}_A)$ equals the
projection of the index value ${\rm Ind}(E^{**}) \in {\rm Z}(A^{**})$
to the discrete part of $A^{**}$, i.~e.~${\rm Ind}(E)$. Consequently,
$F_1({\rm id}_A)$ belongs to $A$ if and only if ${\rm Ind}(E)$ is contained
in the standard injective $*$-representation of $A$ in the discrete part of
$A^{**}$, and $F_1({\rm id}_A)$ equals ${\rm Ind}(E)$.

\smallskip \noindent
If ${\rm Ind}(E)$ is contained in the center of $B$, then the basic
construction of the Jones' tower can be repeated countably many times
by \cite[Th.~3.5, 3.10]{BDH}, since ${\rm Ind}(E_1)={\rm Ind}(E)$ and the
index value is stabilized throughout the tower.
$\: \bullet$

\medskip \noindent
Let us give an example showing that the index value can be outside $A$, and
that ${\rm Ind}(E)$ is very different from ${\rm Ind}(E^{**})$, in general.

\begin{example}  {\rm
  Let $A={\rm C}(S^1)$ be the C*-algebra of all continuous functions on the
  unit circle $S^1= \{ {\rm e}^{{\rm i}\phi} : \phi \in [0,2\pi)\}$, where
  $S^1$ is equipped with the usual topology. Consider the conditional
  expectation
  \[
  E(f)(x) = \frac{f(x)+f(\overline{x})}{2} \, , \: (x \in S^1) \, ,
  \]
  for $f \in A$, where $\overline{x}$ denotes the complex conjugate of $x \in
  S^1$. Obviously, the mapping $(2 \cdot E - {\rm id}_A)$ is positive on $A$
  and $L(E)=K(E)=2$. There does not exist any finite quasi-basis in the sense
  of Y.~Watatani. The discrete part of the bidual linear space and W*-algebra
  $A^{**}$ of $A$ is isomorphic to $l_\infty(S^1)$ by the Gel'fand theorem.
  The formula defining $E^{**}$ is the same as for $E$. The index can be
  counted by \cite[Th.~3.5]{BDH}. The value is
  \[
  {\rm Ind}(E)= \left\{ \begin{array}{ccc}
                        2 & : & x^2 \not= 1 \\ 1 & : & x^2=1
                        \end{array} \right.  \in l_\infty(S^1) \, ,
  \]
  and ${\rm Ind}(E) \not\in A$. Since $A^{**}$ is commutative (hence, type I)
  let us have a look on the non-discrete part of it. By the decomposition
  theory of commutative W*-algebras it is a direct sum of W*-algebras
  $L^\infty([0,1],\lambda)$, where $\lambda$ denotes the Lebesgue measure,
  cf.~\cite{Take}. Consider the canonical embedding $A={\rm C}(S^1) \subset
  L^\infty(S^1,\lambda) \cong L^\infty([0,1],\lambda)$. Again, the index
  of $E^{**}$ reduced to $L^\infty(S^1,\lambda)$ can be found, it equals
  $f(x) \equiv 2$, which is different from the value obtained for the discrete
  part of ${\rm Ind}(E^{**})$.    }
\end{example}

\noindent
The previous example is closely related to a Stinespring theorem for
conditional expectations. We found a more general result in a recent
preprint of G.~J.~Murphy \cite[Th.~2.4]{Mu2}, but the construction is rather
different and more complicated because of the greater ge\-nerality grasped
there.

\begin{theorem}
  Let $E: A \to B \subseteq A$ be a conditional expectation of finite index.
  Then there exists a Hilbert $B$-module $\{ {\cal M}, \langle .,. \rangle \}$,
  a $*$-homomorphism $\pi: A \to {\rm End}_B({\cal M})$ and a partial isometry
  $V$ mapping the Hilbert $B$-module $\{ A, E(\langle .,. \rangle_A) \}$ to
  $\cal M$ so that $E(a)=$ $=V^* \circ \pi(a) \circ V$ for every $a \in A$  and
  the linear hull of the set $\{ \pi(a_1)(V(a_2)) : a_1,a_2 \in A \}$ is
  norm-dense in $\cal M$.
\end{theorem}

\noindent
{\it Proof:} First of all we have to construct a suitable right Hilbert
$B$-module $\cal M$ to represent $A$ on it. We consider $A$ as a right
$B$-module multiplying with elements of $B$ from the right. Start with the
algebraic tensor product $A \odot A$ and factorize the resulting linear
space by the kernel of the $B$-valued inner pre-product
\[
\langle a \otimes x , b \otimes y \rangle = E(x^* E(a^*b)y) \, , \;
a,b,x,y \in A \, .
\]
The completion of the resulting linear factor space with respect to the norm
derived from this $B$-valued inner product is denoted by $\cal M$.
Define the $*$-homomorphism $\pi: A \to {\rm End}_B({\cal M})$ by the formula
$\pi(a)(x \otimes y) = ax \otimes y$ for $a,x,y \in A$. It is a faithful
representation of $A$. The partial isometry $V$ mapping the Hilbert $B$-module
$\{ A, E(\langle .,. \rangle_A) \}$ to $\cal M$ can be obtained by the rule
$V(x) = 1_A \otimes x$ for $x \in A$. Respectively, $V^*(x \otimes y) = E(x)y$
for $x,y \in A$. Consequently, for each $a \in A$ the element $V^* \circ
\pi(a) \circ V$ acts as the $B$-linear operator $E(a)$ on $A$ 
multiplying elements of $A$ by $E(a)$ from the left.
$\: \bullet$

\section{Further results}

\noindent
We would like to describe the local action of conditional expectations of
finite index in a similar way as we have been able to do it in the discrete
W*-case in section one. But, we cannot give such a description on the
elements of the C*-algebra $A$ where the conditional expectation $E$ acts, in
general. Using a fundamental principle of non-commutative topology we should
switch from the minimal projections of the bidual Banach space and
W*-algebra $A^{**}$ either to the maximal (modular) one-sided ideals of $A$
or to the maximal (modular) hereditary C*-subalgebras of $A$ or to the
pure states of $A$.

\begin{corollary}
Let $A$ be a C*-algebra and $E$ be a conditional expectation
on it for which there exists a number $K \geq 1$ such that the mapping
$(K \cdot E - {\rm id}_A)$ is positive. Then the preimage of every maximal
(modular, norm-closed) left ideal of the image C*-algebra $B=E(A) \subseteq A$
contains the set union of at most $[K(E)]$ maximal (modular, norm-closed) left
ideals of $A$ with pairwise orthogonal complements of their w*-carrier
projections, where $[K(E)]$ denotes the integer part of $K(E)$. The
intersection of these maximal left ideals of $A$ with $B$ gives that maximal
left ideal of $B$ back we started with.
\end{corollary}

\noindent
{\it Proof:}  Let $p \in B^{**}$ be minimal. Recall, that the carrier
projections of maximal (modular, norm-closed) one-sided ideals of the
C*-algebra $A$ (which are contained in  the bidual W*-algebra $A^{**}$)
correspond to the minimal projections of the bidual
W*-algebra $A^{**}$ one-to-one by taking orthogonal complements.
By Corollary \ref{cor99},(i) $p$ is a finite sum of minimal projections
$\{q_\alpha \} \in A$, which corresponds to the relation of w*-closed left
ideals
\[
\cap_\alpha q_\alpha A^{**} \equiv p A^{**}
\]
for every such decomposition. Moreover, $E(p A^{**}) = p B^{**}$ and
$E(q_\alpha A^{**}) \subseteq p B^{**}$ for every index $\alpha$. Suppose
$E(r A^{**}) \subseteq p B^{**}$ for a minimal projection $r \in A^{**}$.
Then $E(r) = \mu p$ with $\mu \in (0,1]$ since $p \in B^{**}$ is minimal and
$E$ is faithful. This implies $r \leq p$, and $r \in pA^{**}p$, which is a
finite dimensional C*-algebra of dimension not greater than $[K(E)]^2$ by
Corollary \ref{cor99},(i). Hence, the number of minimal w*-closed left
ideals of $A^{**}$ mapped into the minimal w*-closed ideal $p B^{**}$ of
$B^{**}$ which possess pairwise orthogonal carrier projections
is limited by the number $[K(E)]^2$, as well as the number of pairwise
orthogonal minimal projections $\{ q_\alpha \} \in A^{**}$ mapped to a
multiple of $p \in B^{**}$.
Applying again Akemann's correspondence and observing that $E$ preserves the
maximality and the ideal property of maximal one-sided ideals of
$A$ we are done.
$\: \bullet$

\smallskip \noindent
Because of the interrelation of pure states, maximal (modular, norm-closed)
left ideals and maximal (modular) hereditary C*-subalgebras for
arbitrary C*-algebras we derive the following fact, cf.~\cite{Pe,Mu} and
\cite{Grae}:

\begin{corollary}
Let $A$ be a C*-algebra and $E$ be a conditional expectation
on it for which there exists a number $K \geq 1$ such that the mapping
$(K \cdot E - {\rm id}_A)$ is positive. Then the preimage of every maximal
(modular) hereditary C*-subalgebra of $B=E(A) \subseteq A$ contains
the set union of at most $[K(E)]$ maximal (modular) hereditary C*-subalgebras
$p_\alpha A^{**} p_\alpha \cap A$ of $A$ with pairwise orthogonal complements
of their w*-carrier projections $p_\alpha$. \newline
Every pure state on $B$ has at most $[K(E)]$ different extensions to pure
states on $A$ such that their restrictions to $B \subseteq A$ equal the
original pure state and that they possess pairwise orthogonal w*-carrier
projections.
\end{corollary}


\medskip \noindent
As a by-product of our considerations we reobtain a fact which was observed
in the case of type ${\rm II}_1$ factors by V.~F.~R.~Jones \cite[p.~6]{Jones}
and which was proven for the general W*-case by E.~Andruchow and D.~Stojanoff
\cite[Cor.~2.4]{AS2} and S.~Popa \cite{Po}. Moreover, we can estimate the
dimension of the relative commutant $N' \cap M$ in terms of the constants
$K(E)$ and ${\rm dim}({\rm Z}(N))$.

\begin{corollary}
  Let $E:M \rightarrow N$ be a normal conditional expectation on a W*-algebra
  $M$ such that there exists a number $K \geq 1$ for which the mapping
  $(K \cdot E - {\rm id}_A)$ is positive. Then the center of $M$ is
  finite-dimensional if and only if the center of $N$ is finite-dimensional if
  and only if the relative commutant $N' \cap M$ is finite-dimensional,
  and
   \[
   {\rm dim}(N' \cap M) \leq [K(E)]^2 \cdot {\rm dim}({\rm Z}(N)) \, .
   \]
  (Note, ${\rm dim}({\rm Z}(N)) \leq {\rm dim}(N' \cap M)$, ${\rm dim}
  ({\rm Z}(M)) \leq {\rm dim}(N' \cap M)$ because ${\rm Z}(N) \subseteq (N'
  \cap M)$, ${\rm Z}(M) \subseteq (N' \cap M)$.)
\end{corollary}

\noindent
{\it Proof:} The equivalence of the three claimed conditions follows from
\cite[Cor.~3.19]{BDH} and from Theorem \ref{theorem} immediately.
Since the center of the relative commutant $N' \cap M$ 
contains ${\rm Z}(N)$ and since the minimal projections of ${\rm Z}(N)$
commute with $E$ the conditional expectation $E$ is the direct sum of
${\rm dim}({\rm Z}(N))$ states on $N' \cap M$ with non-intersecting areas
of definition. But, for states on matrix algebras and conditional expectations
of finite index on them the dimension of the matrix algebra is bounded by
$[K(E)]^2$ for $K(E)$ being the structural constant of $E$, cf.~Examples
\ref{ex1}, \ref{ex2}. This gives the argument.
$\: \bullet$

\medskip \noindent
Our next result generalizes \cite[Cor.~2.3]{AS2} and gives some estimates
of the dimensions of $M$ and $N$:

\noindent
\begin{corollary}
  Let $M$ be a C*-algebra, $N$ be a finite-dimensional C*-subalgebra of $M$
  and $E:M \rightarrow N \subseteq M$ be a faithful conditional expectation.
  Then $M$ is finite-dimensional if and only if $E$ is of finite index,
  and the estimate
  \[
  {\rm dim}(M) \leq [K(E)]^2 \cdot {\rm dim}(N)^2
  \]
  is valid. If $M$ and $N$ are commutative, then we have the estimate
  \[
  {\rm dim}(M) \leq [K(E)] \cdot {\rm dim}(N) \, .
  \]
\end{corollary}

\noindent
{\it Proof:}
Obviously, if $M$ is finite-dimensional, then it is finitely generated as a
$N$-module and hence, $E$ is of finite Watatani index. To show the
converse let $B$ be a maximal commutative C*-subalgebra of $N$.
Let $A$ be a maximal commutative C*-subalgebra of $M$ containing $B$. Then
$E(A)=B$ since $E(a)b=E(ab)=E(ba)=bE(a)$ for every $b \in B$, every $a \in A$
and $B$ is maximal commutative in $N$ by assumption. Similar to the proof of
Proposition \ref{prop2} every minimal projection of $B$ has not more than
$[K(E)]$ minimal projection-summands in its decomposition inside $A$. This
shows the statement for commutative C*-algebras.

\noindent
By the general theory of C*-algebras the estimate ${\rm dim}(M) \leq
{\rm dim}(A)^2$ is valid for every C*-algebra $M$ and for every
maximal commutative C*-subalgebra $A$ of $M$. 
$\: \bullet$

\medskip \noindent
Applying a result of R.~V.~Kadison \cite{Kad94} we can show the following:

\begin{proposition}
    Let $E: A \to B \subseteq A$ be a conditional expectation of finite index.
    Then the inequality
    \[
    0 \leq (E(a)-a)^2 \leq (K(E)-1) \cdot (E(a^2)-E(a)^2)
    \]
    holds for self-adjoint elements $a \in A$.   Moreover, if
    $p \in A$ is a projection, then $E(p)$ is a projection if and only
    if $E(p)=p \in B$. Otherwise, $E(p)$ has a spectral value
    $0 < \lambda < 1$. Beside this, $E(a^2)=E(a)^2$ if and only if
    $E(a)=a \in B$, and $K(E)=1$ if and only if $E={\rm id}_M$.
    Furthermore, $K(E) \in \{ 1 \} \cup [2,\infty)$.
\end{proposition}

\noindent
{\it Proof:} We apply the theorem on page 29 of \cite{Kad94} to our
situation. R.~V.~Kadison states that for positive mappings $\psi$ between
unital C*-algebras with $\psi(1) \leq 1$ the inequality
\linebreak[4]
$\psi(a)^2 \leq \psi(a^2)$ holds for every self-adjoint element $a$. In our
situation the equality
\linebreak[4]
$(K(E) \cdot E-{\rm id}_A)(1_A) = (K(E)-1)1_A$ holds.
Dividing by $(K(E)-1)$ we obtain the inequality
\[
\left( (K(E) \cdot E-{\rm id}_A)(a) \right)^2 \leq
(K(E)-1) \cdot \left( (K(E) \cdot E(a^2)-a^2 \right)
\]
which is valid for every self-adjoint $a \in A$. After some obvious
transformations we arrive at
\[
K(E)^2 \cdot E(a)^2-K(E) \cdot aE(a) -K(E) \cdot E(a)a \leq
K(E) \cdot \left( (K(E)-1) \cdot E(a^2)-a^2 \right) \, ,
\]
an inequality which can be transformed to the inequality claimed above.
Setting $E(p)=E(p)^2$ for a projection $p \in A$ we derive $0=(E(p)-p)^2$
and hence, $E(p)=p$. Since in general $1_A \geq E(p)=E(p^2) \geq E(p)^2 \geq
0$ for every projection $p \in A$ the condition $p \not\in B$ implies
the existence of spectral values of $E(p)$ strongly between zero and one.

\noindent
The case $K(E)=1$ is obvious. Applying $E$ to the derived above inequality
and investigating the resulting inequality for $K(E) \in [1,2)$ we obtain
$K(E) \in \{ 1 \}$ as the only possibility. Example \ref{ex1} realizes any
possible value of $K(E)$ inside $[2,\infty)$.
$\: \bullet$

\bigskip \noindent
{\bf Acknowledgement:} The authors are grateful to E.~Andruchow, Y.~Denizeau
and J.-F.~Havet, E.~Christensen, G.~A.~Elliott, R.~Schaflitzel, E.~Scholz and
{\c{S}}.~Str{\u{a}}til{\u{a}} for helpful comments and discussions on the
subject.


\begin{tabular}{lp{0cm}l}
Fed.~Rep.~Germany  & &        Fed.~Rep.~Germany\\
Universit\"at Leipzig & & Humboldt-Universit\"at zu Berlin\\
FB Mathematik/Informatik & & FB Mathematik\\
Mathematisches Institut & & Institut f\"ur reine Mathematik\\
Augustusplatz 10 & & Ziegelstr.~13 A \\
D-04109 Leipzig & & D-10117 Berlin \\
frank@mathematik.uni-leipzig.de & & kirchbrg@mathematik.hu-berlin.de\\
 & & \\
in 1998 in leave to: & & \\
Department of Mathematics & & \\
University of Houston & & \\
Houston, TX 77204-3476, U.S.A. & & \\
frank@math.uh.edu & & 
\end{tabular}

\end{document}